\documentclass[11pt]{article}
\usepackage{epic}
\usepackage{graphics}
\usepackage{graphicx}
\usepackage{eepic}
\usepackage{amsfonts}
\begin{document}
%\setbox1\vbox{\begin{center}
%%\vskip-5mm
%\small\hskip -1cm D\'epartement de Math\'ematiques\\ \vskip-2pt
%\small \hskip -1cm  MAPMO\\ \vskip-2pt
%\small\hskip -1cm  Universit\'e d'Orl\'eans\\ \vskip-2pt
%\small \hskip -1cm  BP 6759\\ \vskip-2pt
%\small\hskip -1cm   F-45067 Orl\'eans cedex 2\\ \vskip-1.5pt
%\small\hskip -1cm   e-mail: Athanasios.Batakis@labomath.univ-orleans.fr
%\end{center}}
%\ht1=10mm \dp1=-5mm \wd1=146mm
\title{Dimension of the harmonic measure of non-homogeneous Cantor sets}
\author{\large\bf Athanasios  B\normalsize\bf{ATAKIS}} 
% {\box1}
\date{}
\maketitle
\newtheorem{Theoreme}{Th\'eor\`eme}[section]
\newtheorem{Theorem}{Theorem}[section]
\newtheorem{Th}{Th\'eor\`eme}[section]
\newtheorem{Def}[Th]{D\'efinition}
\newtheorem{Defi}[Theorem]{Definition}
\newtheorem{Pro}[Th]{Proposition}
\newtheorem{Lemma}[Theorem]{Lemma}
\newtheorem{Corollary}[Theorem]{Corollary}
\newtheorem{Proposition}[Theoreme]{Proposition}
\newtheorem{Lemme}[Theoreme]{Lemme}
\newtheorem{Corollaire}[Theoreme]{Corollaire}
\newtheorem{Consequence}[Theoreme]{Cons\'equence}
\newtheorem{Remarque1}[Theoreme]{Remarque}
\newtheorem{Convention}[Theoreme]{{\sc Convention}}
\newtheorem{PP}[Theoreme]{Propri\'et\'es}
\newtheorem{Conclusion}[Theoreme]{Conclusion}
\newtheorem{Ex}[Theoreme]{Exemple}
\newtheorem{Definition}[Theoreme]{D\'efinition}
\newtheorem{Remark1}[Theorem]{Remark}
\newtheorem{Not}[Theoreme]{Notation}
\newtheorem{Nota}[Theorem]{Notation}
\newtheorem{Propo}[Theorem]{Proposition}
\newtheorem{Exe1}[Th]{Exemple}
\newtheorem{PPtes}[Th]{Propri\'et\'es}
\newenvironment{Proof}{\medbreak{\noindent\bf Proof }}{~{\hskip3pt$\bullet$\bigbreak}}\newenvironment{Proof1}{\medbreak{\noindent\bf Sketch of Proof }}{~{\hskip3pt$\bullet$\bigbreak}}
\newenvironment{Demonstration}{\medbreak{\noindent\bf D\'emonstration
    }}{~{\hskip 3pt$\bullet$\bigbreak}}
\newenvironment{Remarque}{\begin{Remarque1}\em}{\end{Remarque1}}
\newenvironment{Remark}{\begin{Remark1}\em}{\end{Remark1}}
\newenvironment{Remarques}{\begin{Remarque1}\em \ \\* }{\end{Remarque1}}\renewcommand{\Re}{{\cal R}}
\renewcommand{\Im}{{\cal F}}
\newcommand{\finpreuve}{~{\hskip 3pt$\bullet$\bigbreak}}
\newcommand{\hp}{\hskip 3pt}
\newcommand{\hph}{\hskip 8pt}
\newcommand{\hphh}{\hskip 15pt}
\newcommand{\vp}{\vskip 3pt}
\newcommand{\vpv}{\vskip 15pt}
\newcommand{\Omo}{\Omega}
\newcommand{\dist}{\mbox{dist}}
\newcommand{\rd}{{\mathbb R}^2}
\newcommand{\Reel}{{\mathbb R}}
\newcommand{\R}{{\mathbb R}}
\newcommand{\reel}{{\mathbb R}}
\newcommand{\Hyper}{{\mathbb H}}
\newcommand{\Int}{{\mathbb I}}
\newcommand{\Kgros}{I\hskip -4pt K}
\newcommand{\K}{{\mathbb K}}
\newcommand{\Boule}{{\mathbb B}}
\newcommand{\N}{{\mathbb N}}
\newcommand{\Esp}{{\mathbb E}}
\newcommand{\Complex}{{\mathbb C}}
\newcommand{\Ha}{{\mathcal H}}
\newcommand{\Pa}{{\mathcal T}}
\newcommand{\Harm}{{\bold H}}
\newcommand{\Lcal}{{\cal L}}
\newcommand{\ds}{\displaystyle}
\newcommand{\un}{{\mathbf 1}}
\newcommand{\Cone}{C(x,r,\epsilon,\theta)}
\newcommand{\Cn}{C(x,2^{-n},\epsilon ,\theta)}
\newcommand{\Tranche}{W(x,r,\epsilon,\theta)}
\newcommand{\Wn}{W(x,2^{-n},\epsilon,\theta)}
\newcommand{\WFn}{W(x,2^{-n},\epsilon,\theta)\cap F}
\newcommand{\ovec}{\overrightarrow}
\newcommand{\red}{{\bold R}}
\newcommand{\dimH}{\dim_{\Ha}}
\newcommand{\Dim}{\mbox{Dim}}
\newcommand{\diam}{\mbox{diam}}
\newcommand{\diamit}{\mbox{\em diam}}
\newcommand{\para}{\vskip 2mm}
\newcommand{\cod}{\stackrel{\mbox{\tiny cod}}{\sim}}
\newcommand{\cardit}{\mbox{\em card}}
\newcommand{\card}{\mbox{card}}
\newcommand{\Sphere}{{\mathbb S}_d}
\newcommand{\distit}{\mbox{\em dist}}
\begin{abstract}
We prove that the dimension of the
harmonic measure of the complementary of a translation-invariant type of Cantor sets as a continuous function of the
parameters determining these sets. This results extend a previous one of the author and do not use ergodic theoretic tools, not applicable to our case.
\end{abstract}

\section{Introduction}
The purpose of this work is to complement the study of the dimension of the harmonic
measure of the complementary of (not necessarily
self-similar) Cantor sets as a function of  parameters
assigned to these sets. In a previous work \cite{Bata2} we have proved that the parameters assigned to self-similar Cantor sets are continuity points for this function. A new method allows us to treat the continuity over the entire family of parameters determining these translation invariant Cantor sets. We restrain ourselves to sets in the plane for convenience, even
though the proof can be applied to all ``translation-invariant'' Cantor sets in
$\Reel^n$, $n\ge 2$.

Let us start by recalling the definition of the Hausdorff dimension of a measure;
we will use the notation $\dim_{\Ha}$ for the Hausdorff dimension of sets.
\begin{Defi}
If $\mu$ is a measure on $\K$, we will denote by $\dim_*(\mu)$ the lower Hausdorff dimension of $\mu$:
$$\dim_*{\mu}=\inf\{\dim_{\Ha} E\; ;\; E\subset\K\mbox{ and }\mu(E)>0\}$$ and by 
$\dim^*(\mu)$ the upper Hausdorff dimension of $\mu$:
$$\dim^*{\mu}=\inf\{\dim_{\Ha} E\; ;\; E\subset\K\mbox{ and }\mu(\K\setminus E)=0\}.$$
If, for a measure $\mu$ on $\K$, we have $\dim^*(\mu)=\dim^*(\mu)$ then we note this common value $\dim(\mu)$. In the latter case the measure is called exact.
\end{Defi}

For convenience and in order to fix ideas we consider a particular case of translation invariant Cantor sets; we study $4$-corner Cantor sets  constructed in the
following way (see also \cite{Batakis}): let $\underline A,\overline A$ be two constants with
$0<\underline A\le\overline A <\frac{1}{2}$  and let
$(a_n)_{n\in\N}$ be a sequence of real numbers 
with $\underline A\le a_n\le \overline A $ for all $n\in\N$. 

\begin{figure}
\centerline{\setlength{\unitlength}{0.0140in}
\begingroup\makeatletter\ifx\SetFigFont\undefined
% extract first six characters in \fmtname
\def\x#1#2#3#4#5#6#7\relax{\def\x{#1#2#3#4#5#6}}%
\expandafter\x\fmtname xxxxxx\relax \def\y{splain}%
\ifx\x\y   % LaTeX or SliTeX?
\gdef\SetFigFont#1#2#3{%
  \ifnum #1<17\tiny\else \ifnum #1<20\small\else
  \ifnum #1<24\normalsize\else \ifnum #1<29\large\else
  \ifnum #1<34\Large\else \ifnum #1<41\LARGE\else
     \huge\fi\fi\fi\fi\fi\fi
  \csname #3\endcsname}%
\else
\gdef\SetFigFont#1#2#3{\begingroup
  \count@#1\relax \ifnum 25<\count@\count@25\fi
  \def\x{\endgroup\@setsize\SetFigFont{#2pt}}%
  \expandafter\x
    \csname \romannumeral\the\count@ pt\expandafter\endcsname
    \csname @\romannumeral\the\count@ pt\endcsname
  \csname #3\endcsname}%
\fi
\fi\endgroup
\begin{picture}(195,209)(0,-10)
\put(159,30){4}
\put(180,6){4}
\put(132,4){3}
\put(178,53){2}
\path(70,172)(70,194)(47,194)
	(47,172)(70,172)
\path(70,125)(70,147)(47,147)
	(47,125)(70,125)
\path(22,125)(22,147)(0,147)
	(0,125)(22,125)
\path(148,172)(148,194)(125,194)
	(125,172)(148,172)
\path(195,172)(195,194)(173,194)
	(173,172)(195,172)
\path(195,125)(195,147)(173,147)
	(173,125)(195,125)
\path(148,125)(148,147)(125,147)
	(125,125)(148,125)
\path(148,47)(148,69)(125,69)
	(125,47)(148,47)
\path(195,47)(195,69)(173,69)
	(173,47)(195,47)
\path(195,0)(195,22)(173,22)
	(173,0)(195,0)
\path(148,0)(148,22)(125,22)
	(125,0)(148,0)
\path(70,0)(70,22)(47,22)
	(47,0)(70,0)
\path(70,47)(70,69)(47,69)
	(47,47)(70,47)
\path(22,47)(22,69)(0,69)
	(0,47)(22,47)
\path(22,0)(22,22)(0,22)
	(0,0)(22,0)
\path(22,172)(22,194)(0,194)
	(0,172)(22,172)
\put(8,180){1}
\put(132,130){3}
\put(52,179){2}
\put(54,133){4}
\put(7,130){3}
\put(32,157){1}
\put(157,158){2}
\put(134,180){1}
\put(180,179){2}
\put(178,130){4}
\put(54,54){2}
\put(31,30){3}
\put(8,55){1}
\put(7,7){3}
\put(54,7){4}
\put(132,53){1}
\path(195,0)(195,194)(0,194)
	(0,0)(195,0)
\path(70,0)(70,69)(0,69)
	(0,0)(70,0)
\path(195,0)(195,69)(125,69)
	(125,0)(195,0)
\path(195,125)(195,194)(125,194)
	(125,125)(195,125)
\path(70,125)(70,194)(0,194)
	(0,125)(70,125)
\end{picture}}
\caption{A $4$-corner Cantor set and its enumeration}\label{Cantor.fig}
\end{figure}
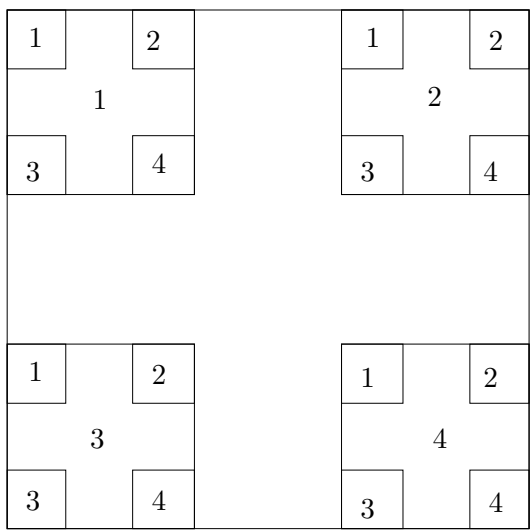

We replace  the square $[0,1]^2$ by four squares of sidelength $a_1$ 
situated in the four corners of $[0,1]^2$. Each of these
squares is then replaced by four squares of sidelength $a_1a_2$
situated in its four corners. At the $n$th stage of the construction 
every square of the $(n-1)$th generation will be replaced by four
squares of sidelength 
$a_1...a_n$ situated in its four corners (see figure
\ref{Cantor.fig}). Let $\K$ be the Cantor set constructed by
repeating the procedure. 

Recall that the harmonic measure of a domain is supported by its boundary and can be seen as the distribution of the exit 
points of Brownian motion starting at some (any) point of the domain (for more details see \cite{Doob}, \cite{Helms} and \cite{Brelot}). 
Carleson  \cite{Ca} has shown that for self-similar $4$-corner Cantor
sets (the sequence $(a_n)_{n\in\N}$ is constant)  
the dimension of the harmonic measure of their complementary 
is strictly smaller than $1$. His proof, involving ergodic theory
techniques, was improved by  Makarov and Volberg
\cite{MV} who showed that the dimension of the
harmonic measure of any self-similar $4$-corner  Cantor set is strictly
smaller than the dimension of the Cantor set. Volberg (\cite{V1}, \cite{V2}) extended these results to a class of
dynamic Cantor repellers. Other comparisons of harmonic and maximal
measures for dynamical systems  are proposed in \cite{BPV}, \cite{LyB}, \cite{UZ}. 
More recently, a multifractal study of harmonic measure on simply connected domains and on Julia sets of polynomial mappings is carried out in \cite{makarov2}, \cite{BMS}.  

In \cite{Batakis}
it is shown that the  dimension of the harmonic measure of the
complementary of $4$-corner Cantor  sets is strictly smaller than
the Hausdorff dimension of the Cantor set, even when the sequence
$(a_n)_{n\in\N}$ is
not constant. In \cite{These} we prove that small
perturbations of the sidelength of the squares of the construction of
$\K$ do not alterate this property.
This last result can also be seen as an immediate consequence of the following theorem.

\begin{Theorem}\label{Th1}
Let $\K=\K_{(a_n)}$  be the $4$-corners Cantor set associated to a 
sequence $a_n$ and  $\K'=\K_{(a_n')}$ a second Cantor set of the same
type associated to the sequence
$(a_n')_{n\in\N}$. Let $\omega$ 
and $\omega'$ be the harmonic measures of $\R^2\setminus\K$ and
$\R^2\setminus\K'$ respectively. Then for all $\epsilon>0$ there
exists a
$\delta=\delta(\epsilon,\overline{A},\underline{A})>0$ such that if  $|a_n'-a_n|<\delta$ for all
$n\in\N$ then 
$|\dim\omega-\dim\omega'|<\epsilon$.
\end{Theorem}
When the sequence $(a_n)_{n\in\N}$ is constant the partial result is already established in \cite{Bata2} using ergodic theoretic tools, which are not applicable in the general case.
\begin{Remark}\label{remark}
Let ${\cal D}:\ell^{\infty}([\overline{A},\underline{A}])\to [0,1]$ be the function that assigns to a sequence  $(a_n)_{n\in\N}\subset [\overline{A},\underline{A}]$ the dimension of harmonic measure of the Cantor set associated to $(a_n)_{n\in\N}$. By refining the estimations in the demonstration of the theorem, we can even show that ${\cal D}$ is a Lipschitz continuous function. The proof of this statement is very technical but straightforward and therefore ommited. 

In particular this refinement implies that if $\ds\sum_{n\in\N}|a_n-a_n'|<\infty$ then the harmonique measures of the corresponding Cantor sets are of the same dimension, and even equivalent when we report them onto the abstract Cantor set $\{1,2,3,4\}^{\N}$.
\end{Remark}
\section{Notations and Preliminary results}
In this section we  establish some estimates on the harmonic
measure of a Cantor set under perturbation, and recall some known
results on the harmonic measure of Cantor-type sets. We also introduce
the tools needed, such as the Hausdorff dimension and the entropy of a
probability measure on a Cantor set.

Let $\K$ be a $4$-corner  Cantor set as described in the
introduction. 
We enumerate $\K$ by identifying it to the abstract Cantor set
$\{1,...,4\}^{\N}$. 
We denote 
$  I _{i_1...i_n}$, where $i_j\in\{{1,2,3,4\}}$ for $ 1\le j\le
n$, the $4^n$ squares  
of the $n$-th generation of the construction of $\K$ with the
enumeration 
shown in the figure \ref{Cantor.fig} and 
the usual condition that 
$  I _{i_1...i_n}$ is the ``father" of the sets $  I
_{i_1...i_ni}$ ,  
$i\in \{{1,2,3,4\}}$. It is clear that  
${\ds \overline A\ge{{\mbox{diam}  I _{i_1...i_ni}}\over 
{\mbox{diam}  I_{i_1...i_n}}}=a_{n+1}\ge \underline  A},\hskip
5pt  i=1,...,4$. 

The collection of the squares of the
$n$-th generation of the constuction of $\K$ will be ${\cal
F}_n=\{I_{i_1...i_n}\, ;\, 
i_1,...,i_n=1,...,4\}$, for $n\in\N$.   
For a square $I\in {\cal F}_n$  we note $\widehat I$ the ``father'' of $I$, i.e. the unique square of $\Im_{n-1}$ containing $I$. If 
$I=I_{i_1...i_k}\in\Im_k$ and $J=I_{j_1...j_n}\in\Im_n$  we will
note $IJ=I_{i_1...i_kj_1...j_n}\in 
\Im_{n+k}$. Finally, for $x\in\K$ and $n\in\N$ let $I_n(x)$ be
the unique square of $\Im_n$ containing $x$.

For a domain $\Omega$, a point $x\in\Omega$ and a Borel set
$F\subset\Reel^2$  we denote by $\omega(x,F,\Omega)$ the
harmonic measure of $F\cap\partial\Omega$\ (for the domain $\Omega$)
assigned to the point $x$. 
Clearly, $F$ carries no measure if it does not intersect
$\partial\Omega$. If 
$\Omega$ is not specified  it will be  $\rd\setminus\K$ and if 
$x$ is the  point at infinity we will simply note $\omega(F)$. 
Finally, for a Borel set $E\subset\Reel^2$ we note $\dim E$ the
Hausdorff dimension of the set $E$. 
%Hence, for $i_1,...,i_n\in\{1,...,4\}$ we denote
%$I_{i_1...i_n}$ the square that corresponds to the cylinder
%$$\{i_1,...,i_n,j_{n+1},... \mbox{ with } j_{n+1},...\in \{1,...,4\}\}$$

\subsection{Dimension of measures}
In this section we  recall some known results on the dimensions of
measures (see also \cite{Mattila}, \cite{Bil}, \cite{Fan}, \cite{You}).
One can prove  (see for instance \cite{Fan}, \cite{Heurteaux}, \cite{BH}) that if
$\mu$ is exact, i.e. if $\dim_*\mu=\dim^*\mu$, then 
\begin{equation}
 \ds\dim\mu=\liminf_{r\to 0}\frac{\log\mu B(x,r)}{\log r}\; ,\;
\mu\mbox{-almost everywhere.} 
\end{equation}
If the  probability measure $\mu$ is supported by a $4$-corner Cantor
set, the balls $B(x,r)$ can be replaced by the 
squares  of the construction of the Cantor set (see \cite{Be}, \cite{These}):
\begin{equation}\label{Remarqueprincipale}
 \ds\dim\mu=\liminf_{n\to \infty}\frac{\log\mu
   \bigl(I_n(x)\bigr)}{\log l\bigl(I_n(x)\bigr)}\; ,\; 
\mu\mbox{-almost everywhere,} 
\end{equation}
where $l\bigl(I_n(x)\bigr)$ is the sidelength of the square $I_n(x)$
and 
$\underline A^n\le l(I_n)\le \overline A^n$.
\begin{Remark}\label{dimnonmono}
If  $\mu$ is an arbitrary (not  necessarily
monodimensional) probability measure we get
\begin{equation}
\ds \dim_*\mu\leq\liminf_{r\to 0}\frac{\log\mu\bigl(B(x,r)\bigr)}{\log r}\le\dim^*
\mu \mbox{ \ }\mu\mbox{-almost everywhere}. 
\end{equation}
Moreover $\ds\dim^*\mu=\mbox{supess}_{\mu} \liminf_{r\to
0}\frac{\log\mu\bigl(B(x,r)\bigr)}{\log r}$ and $\ds\dim_*\mu=\mbox{infess}_{\mu} \liminf_{r\to
0}\frac{\log\mu\bigl(B(x,r)\bigr)}{\log r}$ . 
\end{Remark}
Some results of the following section are stated without demonstration since they are already proved in \cite{Batakis} and in \cite{Bata2}.

\subsection{Estimating perturbations of the harmonic measure}

Suppose that the $4$-corner  Cantor set $\K$ is associated to the
sequence $(a_n)_{n\in\N}$ and let $\K'$ be another Cantor
set associated to the sequence $(a_n')_{n\in\N}$. Let $({\cal
F}_n)_{n\in\N}$ 
be the collections of squares associated to $\K$ and $({\cal
F}_n')_{n\in\N}$ those associated to $\K'$. 

For $\ds I\in
 {\cal F}_n$ and $\ds I'\in
 {\cal F'}_m$ we will write  $I\cod I'$ if $n=m$ and if $I$ and
$I'$ have the same encoding (with respect to the identification to the
abstract Cantor set $\{1,2,3,4\}^{\N}$). 

Finally, if $I\subset\R^2$ is a square of sidelength $\ell$  and $c$ is a positive number
we note $c\dot I$ the square of sidelength $c\ell$ having the same barycenter as $I$.
If $\omega$ is the harmonic measure of $\rd\setminus \K$ and
$\omega '$ the harmonic measure of $\rd\setminus \K'$ we have established 
the following theorem. 
\begin{Theorem}\label{mainlemma}(cf. \cite{Bata2})
For all $\epsilon>0$ there exists a $\delta=\delta(\epsilon,\overline{A},\underline{A})>0$ such that
\begin{equation}\label{conjecture}
\ds\sup_{n\in\N}|a_n-a_n'|<\delta \Rightarrow
\left|\frac{\omega(I)}{\omega(\widehat I)} :
\frac{\omega'(I')}{\omega'(\widehat{I'})}-1\right|<\epsilon,
\end{equation} 
for all $\ds I\in \bigcup_{n\in\N} {\cal F}_n$ and  $\ds I'\in 
\bigcup_{n\in\N} {\cal F}_n'$ with {\em $I\cod I'$}.
\end{Theorem}
\begin{figure}
\centerline{\includegraphics[width=10cm, height=8cm]{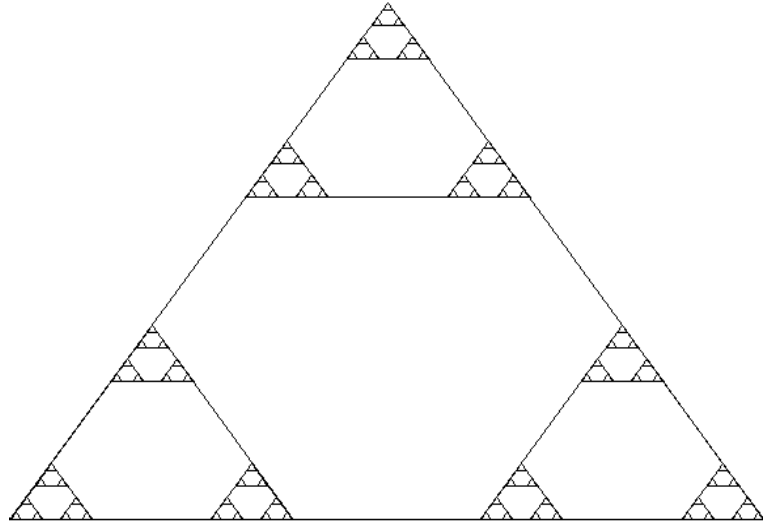}}
\caption{A different type of translation-invariant Cantor set}\label{Cantorgen.fig}
\end{figure}
We will also need the following estimations of the harmonic measure of cylinders (see also \cite{Ca} and \cite{MV} for a version adapted to self-similar sets). The proof uses the ideas already explored in \cite{Bata2}.
\begin{Lemma}\label{HarnackBa}
For every $I,I'\in\Im_n$,
$J\in\Im_k$ and every
$L\in\Im_m$, $n,k,m\in\N$ 
\begin{equation}\label{HarnackCantorauto}
\left|\ds\frac{\omega(IJL)}{\omega(IJ)}:
\frac{\omega(I'JL)}{\omega(I'J)}-1\right|<C\, q^k
\end{equation}
where the constants  $C>0$ and $q\in(0,1)$, depend only on
$\underline{A},\overline{A}$. 
\end{Lemma}
Let us give the proof of this statement.
\begin{Proof}{\bf of lemma \ref{HarnackBa}}
To begin with we need the following Harnack principle (see also \cite{These}, \cite{Ancona2}).
\begin{Lemma}\label{2.5.} (cf. \cite{Ca}, \cite{MV})
Let $\Omega$ be  a domain containing $\infty$ and let 
$A_1 \subset B_1 \subset A_2\subset B_2 \subset ... \subset A_n\subset
B_n$ 
be conformal discs such that the annuli $B_i\setminus A_i$ are
contained in $\Omega$, 
for $1\le i\le n$. If the moduli of the annuli are uniformly bounded
away from zero and if 
$\infty\in  \Omega \setminus B_n$
then, for all pairs of positive harmonic functions 
$u$, $v$ vanishing on  $\partial \Omega \setminus A_1$ and for
all $x \in \Omega \setminus B_n$ we have 
\begin{equation}\label{MVCa}
{\Biggl| {\ds{{u(x)} \over{v(x)}} 
: \ds{{u(\infty)} \over{v(\infty)}}}-1\Biggr|}\leq 
{K  \ds{q^n}}
\end{equation}
where $q<1$ and $K$ are two constants that depend only
on the lower  bound of the moduli of the annuli.
\end{Lemma}
We use this result to prove the following:
%%%%%%%%%%%%%%%%%%%%%%%%%%%%%%%%%%%%%%%%%%%%%%%
\begin{Lemma}\label{Comparaisoncarre2}
There are constants $K>0$ and $0<q<1$ depending only on $\underline A,
\overline A$ such that for all $i,j,k\in\N$ and
for all squares
$I\in\Im_i, J\in\Im_j, K\in\Im_k$ of the 
construction of $\K$, if $Q=c_0 \cdot I$, 
\begin{equation}\label{comparaisoncarre2}
\left|\frac{\omega(x,IJK,Q\setminus \K)}{\omega(x,IJ,Q\setminus \K)} : \frac{\omega(IJK
)}{\omega(IJ)}-1\right|<Kq^j \mbox{
for all } x\in \partial 
\left\{\frac{1+c_0}{2}  \cdot I\right\}.
\end{equation}
The result applies also to the Cantor set $\K'$.
\end{Lemma}
\begin{Proof}
By lemma \ref{2.5.}, 
\begin{equation}\label{RappelHarnack}
\left|\frac{\omega(x,IJK)}{\omega(x,IJ)} :
\frac{\omega(IJK)}{\omega(IJ)}-1\right|<Kq^{j}
\mbox{ , for }   x\notin \frac{1+c_0}{2}\cdot I
\end{equation}
Let $A=\ds\frac{\omega(IJK)}{\omega(IJ)}$. We have
$$\omega(x,IJK,Q\setminus \K)=\omega(x,IJK)-\int_{\partial
Q}\omega(z,IJK)\omega(x,dz,Q\setminus \K),$$
for $x\in  \partial
\left\{\frac{1+c_0}{2}\cdot I\right\}$.

By the equation (\ref{RappelHarnack}), 
$$A\omega(x,IJ)
-Kq^j A \omega(x,IJ)\le \omega(x,IJK)\le
A\omega(x,IJ) +Kq^j A \omega(x,IJ).$$

We get
\begin{eqnarray}
\omega(x,IJK,Q\setminus \K)&\le& A\omega(x,IJ) +Kq^j A
\omega(x,IJ)-\nonumber \\
&-&\int_{\partial Q}\Bigl(A\omega(z,IJ) -Kq^j
A \omega(z,IJ)\Bigr)\omega(x,dz,Q\setminus \K) \nonumber \\
&=& A\omega(x,IJ)-\int_{\partial Q}A\omega(z,IJ)\omega(x,dz,Q\setminus \K)+\nonumber \\
&+ &Kq^j \Bigl(A
\omega(x,IJ) +\int_{\partial Q}A \omega(z,IJ)\omega(x,dz,Q\setminus \K)\Bigr) \nonumber \\
&=&A \omega(x,IJ,Q\setminus \K) +Kq^j \Bigl(A
\omega(x,IJ) +\nonumber\\
&+&\int_{\partial Q}A \omega(z,IJ)\omega(x,dz,Q\setminus \K)\Bigr)
\end{eqnarray}

Therefore, 
\begin{equation}\label{Aterminer}
\ds \frac{\omega(x,IJK,Q\setminus \K)}{\omega(x,IJ,Q\setminus \K)}\le A+Kq^j A\frac{
\omega(x,IJ) +\int_{\partial Q} \omega(z,IJ)\omega(x,dz,Q\setminus \K)}{ \omega(x,IJ,Q\setminus \K)} 
\end{equation}

It suffices now to show that  the quantity 
$$\frac{\omega(x,IJ) +\int_{\partial Q} \omega(z,IJ)\omega(x,dz,Q\setminus \K)}{ \omega(x,IJ,Q\setminus
\K)}$$ is smaller that a given constant.
Take  $x_0 \in\partial
\left\{\frac{1+c_0}{2}\cdot I\right\}$ such that
$$\omega(x_0,IJ)=\max\left\{\omega(x,IJ)\hp ; \hp x\notin
\left\{\frac{1+c_0}{2}\cdot I\right\} \right\}.$$ 

Using the maximum principle we get
\begin{eqnarray*} 
\omega(x_0,IJ,Q\setminus \K)&=&\omega(x_0,IJ)-
\int_{\partial Q}\omega(z,IJ) \omega(x_0,dz,Q\setminus
\K)\\
&\ge& \omega(x_0,IJ)-\int_{\partial Q}\omega(x_0,IJ)
\omega(x_0,dz,Q\setminus \K)\\
&=& \omega(x_0,IJ)\bigl(1-\omega(x_0,\partial Q,Q\setminus
\K)\bigr) 
\end{eqnarray*}
By standard capacitary techniques one can verify  (see \cite{Batakis}) that
$1-\omega(x_0,\partial Q,Q\setminus 
\K)$ is greater that a constant $c>0$ depending only on
$\underline A, \overline A$. 

By using Harnack's principle  we get 
$$1-\omega(x,\partial Q,Q\setminus 
\K)\ge c\mbox{ , for all } x\in\partial
\left\{\frac{1+c_0}{2}\cdot I\right\} ,$$ for a new constant
$c>0$.

Hence, \  $\ds\frac{\omega(x,IJ)+\int_{\partial Q} \omega(z,IJ)\omega(x,dz,Q\setminus \K)}{ \omega(x,IJ,Q\setminus\K)}\le \frac{2}{c}$ \  and therefore, by relation
(\ref{Aterminer}),  
\begin{equation}
\ds \frac{\omega(x,IJK,Q\setminus \K)}{\omega(x,IJ,Q\setminus \K)}\le A(1+\frac{2}{c}Kq^j) 
\end{equation}
On the other hand  $\ds A=\frac{\omega(IJK)}{\omega(IJ)}$; 
% on applying once more  relation  (\ref{RappelHarnack}) 
we obtain
$$ \frac{\omega(x,IJK,Q\setminus \K)}{\omega(x,IJ,Q\setminus \K)} : \frac{\omega(IJK)}{\omega(IJ)}-1<\frac{2}{c}Kq^j\mbox{ , for all } x\in \partial
\left\{\frac{1+c_0}{2}\cdot I\right\},$$ 

The left hand inequality and hence the lemma \ref{Comparaisoncarre2} is established in the same way .
\end{Proof}
It is now evident that 
$$ \frac{\omega(x,IJK,Q\setminus \K)}{\omega(x,IJ,Q\setminus \K)}=\frac{\omega(x,I'JK,Q'\setminus \K)}{\omega(x,I'J,Q'\setminus \K)},$$ for any square $I'\in\Im_n$, where $Q'=c_0\cdot I'$. The proof of lemma \ref{HarnackBa} is complete. 
\end{Proof}
\begin{Corollary}\label{tinvariance}
There is a constant $\tilde C>1$ such that for any $n,k\in\N$, all  $I,I'\in\Im_n$ and every
$J\in\Im_k$ we have
\begin{equation}\label{tinvariancee}
\ds\frac{\omega(IJ)}{\omega(I)}\le \tilde C
\frac{\omega(I'J)}{\omega(I')}
\end{equation}
where the constant  $\tilde C>0$ depends only on
$\underline{A},\overline{A}$. 
\end{Corollary}
The proof of the corollary is an easy application of lemma \ref{HarnackBa}.
%%%Proof
%%%
\section{Proof of the main result.}

This section is dedicated to the proof of theorem \ref{Th1}.
We will  make use of the following known version of the  theorem of
large numbers (see for instance \cite{HH}). 
\begin{Lemma}\label{Grandsnombres}
Let $X_n$ be a sequence of uniformly bounded real random variables on a
probability space 
$({\mathbb X},{\cal B}, P)$ and let $(\Im_n)_{n\in\N}$ be an increasing
sequence of $\sigma$-subalgebra of $\mathbb B$ such that $X_n$ 
is measurable with respect to $\Im_n$ for all
$n\in\N$. Then
\begin{equation}\label{grandsnombres}
\ds \lim_{n\to\infty}\frac{1}{n}\sum_{k=1}^n
\left(X_k-\Esp (X_k|\Im_{k-1})\right)=0\hp  \hp P\mbox{-almost
surely}
\end{equation}
\end{Lemma}
%We give a sketch of proof for completeness. A full  demonstration
%of a more general version is given in \cite{HH}.
%\begin{Proof}
%The random variables $Y_n= X_n-\Esp (X_n|\Im_{n-1})$ are orthogonal
%uniformly bounded by a constant $L>0$ and verify $E(Y_n)=0$. By
%Chebychev's inequality 
%$$P(|\frac{1}{n}\sum_{k=1}^nY_k|>\epsilon)<\frac{L}{n\epsilon^2}.$$
%Then the sequence $\ds \frac{1}{n^2}\sum_{k=1}^{n^2}Y_k$ converges to $0$
%almost surely and a  standard argument
%implies then that $\ds \frac{1}{n}\sum_{k=1}^nY_k$ 
%converges also to 
%$0$ (see also  \cite{HH}). 
%\end{Proof}
The following elementary lemma is also useful; the proof is left to the reader.
\begin{Lemma}\label{elementary}
Let $\alpha_1,...,\alpha_n$ be real numbers such that $\ds\sum_{i=1}^n\alpha_i=0$. Then, for any choise of real values $h_1,...,h_n$, we have $$\begin{array}{rcl}
\left|\sum_{i=1}^n\alpha_i h_i\right|& \le &\max\left\{\ds \sum_{\{i\; ;\;\alpha_i>0\}}\alpha_i\;,\;-\sum_{\{i\; ;\;\alpha_i<0\}}\alpha_i\right\}\left(\ds \max_{1\le i\le n}h_i -\min_{1\le i\le n}h_i\right)\cr
& = &\ds \sum_{\{i\; ;\;\alpha_i>0\}}\alpha_i  \left(\max_{1\le i\le n}h_i -\min_{1\le i\le n}h_i\right).
\end{array}$$
\end{Lemma}

\begin{Proof}{\bf of theorem \ref{Th1}}
\hskip 3mm For $p\in\N$
consider the sequence of $\sigma$-algebras $(\Re_n)_{n\in\N}$ where
$\Re_n$ 
is generated by $\Im_{np}$. 

The hypothesis of lemma
\ref{Grandsnombres} can be easily verified to hold for the sequence
of random variables  $(X_n^p)_{n\in\N}$ given by $\ds X_n^p(x)=\frac{1}{p}\left|\log\left(\frac{ \omega\left(I_{np}(x)\right)}{\omega\left(I_{(n-1)p}(x)\right)}\right)\right|$ and the sequence of
$\sigma$-algebras 
$(\Re_n)_{n\in\N}$.

We get
\begin{equation}\label{helping}
\ds \lim_{n\to\infty}\frac{1}{n}\sum_{k=1}^n\biggl[
X_k^p-\Esp_{\omega}(X_k^p|\Re_{k-1})\biggr]=0\hp \hp
\omega\mbox{-almost everywhere.} 
\end{equation}

On the other hand,  on  $I\in\ds\Re_{n-1}$, $n\in\N$,
$$\Esp_{\omega}(X_n^p|\Re_{n-1})= \frac{1}{p}\sum_{J\in\Im_p}
\frac{\omega(IJ)}{\omega(I)}
\left|\log\left(\frac{\omega(IJ)}{\omega(I)}\right)\right|.$$
We show that this quantity is almost constant on $x\in\K$ if $p$ is taken sufficiently large.

Take $\epsilon>0$, $n\in\N$ and $I\in \Im_n$. For $j,k\in\N$ we have
\begin{eqnarray}\label{starting}
 &&\sum_{J\in\Im_{j+k}}
\frac{\omega(IJ)}{\omega(I)}\log\left(\frac{\omega(IJ)}{\omega(I)}\right)=\sum_{J\in\Im_{j}} \sum_{K\in\Im_k} \frac{\omega(IJK)}{\omega(I)}\log\left(\frac{\omega(IJK)}{\omega(I)}\right)=\cr
 &=&\sum_{J\in\Im_{j}} \sum_{K\in\Im_k} \frac{\omega(IJK)}{\omega(IJ)}\log\left(\frac{\omega(IJK)}{\omega(IJ)}\right)\frac{\omega(IJ)}{\omega(I)}
+ \sum_{J\in\Im_{j}}\frac{\omega(IJ)}{\omega(I)}\log\left(\frac{\omega(IJ)}{\omega(I)}\right)
\end{eqnarray}
For $L\in\Im_n$ and $k,n,j\in\N$ we note  $$h_k(L)=-\frac{1}{k}\sum_{K\in\Im_k} \frac{\omega(LK)}{\omega(L)}\log\left(\frac{\omega(LK)}{\omega(L)}\right).$$ In particular, we put
$$h_k(IJ)=-\frac{1}{k}\sum_{K\in\Im_k} \frac{\omega(IJK)}{\omega(IJ)}\log\left(\frac{\omega(IJK)}{\omega(IJ)}\right) \mbox{ and }
\Delta_j^k(I)= \max_{J\in \Im_j} h_k(IJ)-\min_{J\in \Im_j}h_k(IJ).$$

We will use the following lemma.
\begin{Lemma}\label{technical}
For all $\epsilon>0$, if $j,k\in\N$ are big enough (depending only on $\underline{A}$ and $\overline{A}$) then for all $n\in\N$ and $I\in\Im_n$ we have $\displaystyle \Delta_j^k(I)<\epsilon$.
\end{Lemma}
We first proceed with the proof of this sub-lemma.
\begin{Proof}{\bf of lemma \ref{technical}}
We can rewrite formula (\ref{starting}):
\begin{equation}\label{goingon}
(j+k)h_{j+k}(I)=\sum_{J\in\Im_{j}}k\frac{\omega(IJ)}{\omega(I)} h_k(IJ) +jh_j(I)
\end{equation}

By applying formula (\ref{goingon}) to a cylinder $I=I_1I_2$ with $I_1\in\Im_{i_1}$ and $I_2\in\Im_{i_2}$ we have 
$$(j+k)h_{j+k}(I_1I_2)= \sum_{J\in\Im_{j}}k\frac{\omega(I_1I_2J)}{\omega(I_1I_2)} h_k(I_1I_2J)+jh_j(I_1I_2)$$
Now take $j$ big enough to have $Kq^j<\epsilon$ and afterwards choose $k$ in order to have that $\ds\frac{j}{k+j}<\epsilon$ . 
Remark that, by lemma \ref{HarnackBa}, $h_k(I_1I_2'J)-h_k(I_1I_2J)<2\epsilon$.
We have 
\begin{eqnarray}\label{next}
\Delta_{i_2}^{k+j}(I_1)&=&\max_{I_2\in \Im_{i_2}} h_{k+j}(I_1I_2)-\min_{I_2\in \Im_{i_2}}h_{k+j}(I_1I_2)\cr
&\leq & 5\epsilon+ \max_{I_2,I_2'\in\Im_{i_2}} \left[\sum_{J\in\Im_{j}}\frac{\omega(I_1I_2J)}{\omega(I_1I_2)}  h_k(I_1I_2J)  - \sum_{J\in\Im_{j}}\frac{\omega(I_1I_2'J)}{\omega(I_1I_2')}  h_k(I_1I_2'J) \right]\cr
&\leq &  10\epsilon  +   \max_{I_2,I_2'\in\Im_{i_2}} \sum_{J\in\Im_{j}} \left(\frac{\omega(I_1I_2J)}{\omega(I_1I_2)}-\frac{\omega(I_1I_2'J)}{\omega(I_1I_2')}\right)h_k(I_1I_2J).  
\end{eqnarray}

We can now apply lemma \ref{elementary} to get 
\begin{equation}\label{elemuse}
\Delta_{i_2}^{k+j}(I_1)\le10\epsilon+ \max_{I_2,I_2'\in\Im_{i_2}} \sum_{J\in S_{i_1}^j(I_2,I_2')} \left( \frac{\omega(I_1I_2J)}{\omega(I_1I_2)}-\frac{\omega(I_1I_2'J)}{\omega(I_1I_2')} \right) \Delta_j^k(I_1I_2)
\end{equation}
where ${\mathcal S}_{I_1}^j(I_2,I_2')$ is the set of cylinders $J\in\Im_j$ such that $\ds \frac{\omega(I_1I_2J)}{\omega(I_1I_2)}>\frac{\omega(I_1I_2'J)}{\omega(I_1I_2')}$.

The following lemma is easy to prove:
\begin{Lemma}\label{facile}
There is a constant $0<\zeta<1$ such that for all $i_1,i_2,j\in\N$ , all $I_{1}\in\Im_{i_1}$ and all $I_{2},I_{2}'\in\Im_{i_2}$ we have 
$$ \ds   \sum_{J\in S_{I_1}^j(I_2,I_2')} \left( \frac{\omega(I_1I_2J)}{\omega(I_1I_2)}-\frac{\omega(I_1I_2'J)}{\omega(I_1I_2')} \right) <\zeta$$
where ${\mathcal S}_{I_1}^j(I_2,I_2')$ is the set of cylinders $J\in\Im_j$ such that $\ds \frac{\omega(I_1I_2J)}{\omega(I_1I_2)}>\frac{\omega(I_1I_2'J)}{\omega(I_1I_2')}$.
\end{Lemma}
The proof follows from the translation invariance of $\omega$ (corollary \ref{tinvariance}). 
\begin{Proof}{\bf of  lemma \ref{facile}}
Remark that by lemma \ref{HarnackBa}, if $\displaystyle C=\frac{K}{1-q}$, we have 
$\ds \frac{\omega(IJ)}{\omega(I'J)}\le C$ for all $I,I'\in\Im_n$ , $n\in\N$,  and  all $J\in\bigcup_{n\in\N}\Im_n$.
Hence, 
$$ \sum_{J\in S_{I_1}^j(I_2,I_2')} \left( \frac{\omega(I_1I_2J)}{\omega(I_1I_2)}-\frac{\omega(I_1I_2'J)}{\omega(I_1I_2')} \right) \le \left(1-\frac{1}{C}\right)\sum_{J\in S_{I_1}^j(I_2,I_2')}  \frac{\omega(I_1I_2J)}{\omega(I_1I_2)}\le 1-\frac{1}{C},$$
which is the lemma conclusion for $\ds \zeta=1-\frac{1}{C}$.
\end{Proof}
By applying lemma \ref{facile} to relation (\ref{elemuse}) we conclude that there is $0<\zeta<1$ depending only on $\underline{A}, \overline{A}$ such that
\begin{equation}\label{expo}
\Delta_{i_2}^{k+j}(I_1)\le10\epsilon+\zeta\Delta_j^k(I_1I_2)
\end{equation}
By repeating the same reasoning if we write $k=k_1+k_2$ we can establish the inequalities 
\begin{equation}\label{expobis}
\Delta_j^{k_1+k_2}(I_1I_2)<10\epsilon+ \zeta \Delta_{k_1}^{k_2}(I_1I_2J).
\end{equation}
Hence,
 $$\Delta_{i_2}^{k+j}(I_1)-\frac{10\epsilon}{1-\zeta}\le\zeta\left(\Delta_{k_1}^{k_2}(I_1I_2J) -\frac{10\epsilon}{1-\zeta}\right)\le\zeta^2\left(\Delta_{k_1}^{k_2}(I_1I_2J)-\frac{10\epsilon}{1-\zeta}\right).$$ 
The sequence $\Delta_l^m(I)$ being uniformly bounded we get, by decomposing again $k_2$ and repeating the procedure, that if $k$ is big enough, $\ds \Delta_{i_2}^{k+j}<\frac{20\epsilon}{1-\zeta}$.

The real constant $\zeta$ depending only on $\underline{A}$, $\overline{A}$, the proof is complete.
\end{Proof}

We now apply lemma \ref{Grandsnombres} to an adapted filtration $\Re_n$ : By the previous lemma we can choose $j,k$ such that $\Delta_j^k(I)<\epsilon$. By formula (\ref{goingon}), for this choice of $j$ and $k$ and for all $n\in\N$ there are constants $c_n$ such that
$$\left|\frac{1}{k+j}\Esp\left\{\sum_{\ell=1}^{k+j}X_{n+\ell}^1\Big|\Im_n\right\}-c_n\right|<\epsilon.$$
By lemma \ref{Grandsnombres} and the relation (\ref{helping}) following it we then deduce (for $p=k+j$) 
$$ \left| \liminf_{n\to\infty}\frac{1}{n(k+j)}\sum_{\ell=1}^{n(k+j)} X_{\ell}^1 -\liminf_{n\to\infty}\frac{1}{n}\sum_{\ell=1}^{n} c_{\ell}\right|<\epsilon\;  , \; \omega \mbox{ - a.e. on }\K$$
This implies that $$\left|\liminf_{n\to\infty}\left|\frac{\log\omega(I_{n(k+j)}(x))}{ \log\left(\prod_{i=1}^{n(k+j)}a_i\right)}\right| - \liminf_{n\to\infty}\left|\frac{k+j}{\log\left(\prod_{i=1}^{n(k+j)}a_i\right)}\sum_{\ell=1}^{n} c_{\ell} \right| \right|<\epsilon\;  , \; \omega\mbox{ - a.e. on }\K.$$
On the other hand, once we have fixed $k,j$ we can use lemma \ref{mainlemma} to choose $\delta$ in a way that, for all $n\in\N$,
$$|c_n-c_n'|<\epsilon,\mbox{ and } \frac{1}{n}\left|\log\left(\prod_{i=1}^na_i\right)-\log\left(\prod_{i=1}^na_i'\right)\right|<\epsilon,$$ 
where $c_n'$ is the same sequence associated to the harmonique measure $\omega'$.   
We can finally use relation (\ref{Remarqueprincipale}) to conclude that $|\dim\omega-\dim\omega'|<4\epsilon$.
\end{Proof}

\section{Consequences and remarks}\label{consequences}
 
It is implicitely proved that the harmonic measure of the sets $\K$ studied here satisfy the relationship $\dim_*\omega= h_*^{\K}(\omega)$, where $$h_*^{\K}(\omega)=\liminf_{n\to\infty}\frac{1}{\ds\log\prod_{i=1}^na_i}\sum_{I\in\Im_n}\log\omega(I)\omega(I),$$ and $(a_n)_{n\in\N}$ is the construction sequence associated to $\K$. This fact is a consequence of the space invariance of $\omega$ and is a key factor in the proof of our results.

It is natural to ask whether the relation (\ref{conjecture}) suffices
to conclude that the dimensions of two measures $\omega$ and
$\omega'$ (not necessarily harmonic) are close. This is not the
case. There are counterexemples (see \cite{These}) even when the
measures are doubling on $(\Im_n)_{n\in\N}$ and exact (cf \cite{Batakis4}).
   
Even if the equality between the Hausdorff dimension and the entropy of the
harmonic measure of the complementary of Cantor sets  plays a crucial role in
the proof of theorem \ref{Th1}, we would like to point out that the measures constructed $\mu$ and $\nu$ in the example of  \cite{Batakis4} satisfy the relation
$$ \left|\frac{\mu(I)}{\mu({\widehat I})}:\frac{\mu(I)}{\mu({\widehat I})}-1\right|<\delta,$$
with $\delta$ as small as we want, as well as the the equalities  $h_*(\mu)=\dim\mu$ and $h_*(\nu)=\dim\nu$.
Nevertheless  $\ds |\dim\mu-\dim\nu|\geq\frac{1}{2}$.

To establish the result claimed in remark \ref{remark} we first need a fine precision of inequalities in theorem \ref{mainlemma} and secondly we need to quantify the dependance on $\epsilon$ of the choice of $k$. In fact, we can find  suitable $k$'s that are bounded by  $-C\log(\epsilon)$, where $C$ is a positive constant depending on $\underline{A}$, $\overline{A}$,  which is sufficient in order to prove the claim.

\bibliographystyle{plain}
\bibliography{biblio}
\hfill
\begin{tabular}{c}
{\bf Athanasios BATAKIS}\\
{MAPMO}\\
{Universit\'e d'Orl\'eans}\\
{BP 6759}\\
{45067 Orl\'eans cedex 2}\\
{FRANCE}\\
{email: batakis@labomath.univ-orleans.fr}
\end{tabular}

\end{document}